\documentclass[twocolumn,amsmath,amssymb,aps,secnumarabic,%
    nofootinbib,groupedaddress]{revtex4}
\usepackage{times,mathptmx,amsthm}
\newtheorem{theorem}{Theorem}[section]
\newtheorem{lemma}[theorem]{Lemma}
\theoremstyle{definition}
\theoremstyle{remark}
\newtheorem*{remark}{Remark}

\newcommand{\R}{\mathbb{R}}

\newcommand{\Z}{\mathbb{Z}}
\newcommand{\F}{\mathbb{F}}
\renewcommand{\tilde}{\widetilde}

\newcommand{\va}{\vec{a}}
\newcommand{\vx}{\vec{x}}
\newcommand{\vp}{\vec{p}}
\newcommand{\tO}{\tilde{O}}
\newcommand{\Tr}{\mathrm{Tr}}
\newcommand{\OA}{\mathrm{OA}}

\newcommand{\Eg}{\emph{E.g.}}

\newcommand{\floor}[1]{\lfloor #1 \rfloor}
\newcommand{\bfloor}[1]{\bigl\lfloor #1 \bigr\rfloor}
\newcommand{\ceil}[1]{\lceil #1 \rceil}
\newcommand{\eq}[2]{\begin{equation}\label{#1}#2\end{equation}}

\newcommand{\eatline}{\vspace{-\baselineskip}}

\newcommand{\thm}[1]{Theorem~\ref{#1}}
\renewcommand{\sec}[1]{Section~\ref{#1}}
\newcommand{\lem}[1]{Lemma~\ref{#1}}

\begin{document}
\title{Numerical cubature using error-correcting codes}

\author{Greg Kuperberg}
\email{greg@math.ucdavis.edu}
\thanks{This material is based upon work supported by the National Science
    Foundation under Grant No. 0306681.}
\affiliation{Department of Mathematics, University of
    California, Davis, CA 95616}

\begin{abstract}
\centerline{\textit{\normalsize Dedicated to W{\l}odzimierz and Krystyna
    Kuperberg on the occasion of their 40th anniversary}}
\vspace{\baselineskip}

We present a construction for improving numerical cubature formulas with equal
weights and a convolution structure, in particular equal-weight product
formulas, using linear error-correcting codes.  The construction is most
effective in low degree with extended BCH codes.  Using it, we obtain several
sequences of explicit, positive, interior cubature formulas with good
asymptotics for each fixed degree $t$ as the dimension $n \to \infty$.  Using a
special quadrature formula for the interval \cite{Kuperberg:moments}, we obtain
an equal-weight $t$-cubature formula on the $n$-cube with $O(n^{\floor{t/2}})$
points, which is within a constant of the Stroud lower bound.  We also obtain
$t$-cubature formulas on the $n$-sphere, $n$-ball, and Gaussian $\R^n$ with
$O(n^{t-2})$ points when $t$ is odd.  When $\mu$ is spherically symmetric and
$t=5$, we obtain $O(n^2)$ points.  For each $t \ge 4$, we also obtain explicit,
positive, interior formulas for the $n$-simplex with $O(n^{t-1})$ points; for
$t=3$, we obtain $O(n)$ points.  These constructions asymptotically improve the
non-constructive Tchakaloff bound.

Some related results were recently found independently by Victoir
\cite{Victoir:cubature}, who also noted that the basic construction more
directly uses orthogonal arrays.
\end{abstract}
\maketitle

\section{General results}

Let $\mu$ be a normalized measure on $\R^n$ with finite moments. A
\emph{cubature formula of degree $t$}, or $t$-cubature formula, for $\mu$ is a
set of points $F = \{\vp_a\} \subset \R^n$ and a weight function $\vp_a \mapsto
w_a \in \R$ such that
$$\int P(\vx) d\mu = P(F) = \sum_{a=1}^N w_a P(\vp_a)$$
for polynomials $P$ of degree at most $t$. (If $n=1$, then $F$ is also called a
\emph{quadrature formula}.)  The formula $F$ is \emph{equal-weight} if the
$w_a$ are all equal; \emph{positive} if $w_a > 0$ for all $a$; and otherwise it
is \emph{negative}. Let $X$ be the support of $\mu$.  The formula $F$ is
\emph{interior} if every point $\vp_a$ is in the interior of $X$; it is
\emph{boundary} if every $\vp_a$ is in $X$ and some $\vp_a \in \partial X$; and
otherwise it is \emph{exterior}.  These properties of cubature formulas are
often abbreviated. \Eg, PI means positive and interior and EB means
equal-weight and boundary.  (Exterior formulas are denoted ``O,'' for
outside.)  An equal-weight formula is abbreviated ``E'' and is also called a
(geometric) \emph{$t$-design} or a \emph{Chebyshev-type formula}.

The main use of a cubature formula is to numerically integrate a function $f$
which is approximately a polynomial. In this application, formulas with many
points or non-explicit points are impractical, exterior formulas are
ill-founded if $f$ is only defined on $X$, and formulas with large negative
weights are ill-conditioned on the class of continuous functions \cite[Ch.
1]{Stroud:calc}.  Thus PI formulas with few points are the best kind.

By Tchakaloff's theorem \cite[p. 61]{Stroud:calc}, every measure $\mu$ on
$\R^n$ has a PI $t$-cubature formula with at most $\binom{n+t}{t}$ points, the
same as the dimension of the vector space of relevant polynomials,
$\R[\vx]_{\le t}$.   (If $\partial X$ has non-zero measure, it may only be a PB
formula.)  Tchakaloff's theorem has a short proof, but it is computationally
non-constructive.  Many known formulas with $n$ small, or with $n$ large and $t
\le 2$, are better than the Tchakaloff bound \cite{Stroud:calc,Cools:encyc}.  
But if $n$ is large, $t \ge 3$, and $\mu$ is reasonably natural, most explicit
formulas in the existing literature are either negative, exterior, or have
exponentially many points.

In this article we present a new method to thin equal-weight cubature formulas
with a convolution structure, in particular product formulas for product
measures.  (By \emph{thinning} a formula, we mean removing some of its points
without reducing its cubature degree.) The thinned formulas are efficient in
high dimensions and low degree.  The method also applies to some non-product
measures that are related to product measures, in particular spheres and
simplices with uniform measure. Victoir \cite{Victoir:cubature} independently
obtained the basic construction when $q=2$, together with some other
generalizations not considered by this author.  However, many of our asymptotic
bounds and derived constructions are new.

If $F$ and $G$ are two cubature formulas, we define their
\emph{convolution} $F * G$ to be their sum as sets, $F + G$. The
weight $w_a$ of $\vp_a$ in $F * G$ is given by a product rule:
$$w_a = \sum_{\substack{\vp_a = \vp_b + \vp_c
    \\ \vp_b \in F, \vp_cb \in G}} w_b w_c.$$
Convolution of cubature formulas is related to convolution of measures in two
ways:  First, it is convolution of measures if cubature formulas are interpreted
as atomic measures.  Second, if $F$ is a $t$-cubature formula for $\mu$ and
$G$ is a $t$-cubature formula for $\nu$, then $F * G$ is a
$t$-cubature formula for $\mu * \nu$.  In particular, product formulas and
product measures are convolutions in independent directions.

We also recall some basic facts from coding theory.  For each prime power $q$,
there is a unique finite field $\F_q$ with $q$ elements.  A \emph{linear
error-correcting code} of length $\ell$, dimension $k$, and distance $t$ over
$\F_q$ is a $k$-dimensional vector subspace of $\F_q^\ell$ such that each
non-zero vector has at least $t$ non-zero coordinates.  It is also called an
$[\ell,k,t]_q$ code.  A code $C$ is a \emph{zero-sum code} if the coordinates
of every $\va \in C$ sum to $0$.

\begin{theorem} Let $t$, $n$, and $\ell$ be positive integers, let $q$ be a
prime power, and let $\mu$ be a measure on $\R^n$. For each $1 \le i \le \ell$,
let $F_i$ be an equal-weight formula with $q$ elements such that the
convolution
$$F = F_1 * F_2 * \ldots * F_\ell$$
is a $t$-cubature formula for $\mu$.  Then an $[\ell,k,t+1]_q$
code $C$ yields a thinning $G \subset F$ with $q^{\ell-k}$ points.  In
addition, if each $F_i$ is centrally symmetric, $t$ is odd, and either $q$ is
odd or $C$ is a zero-sum code, then $C$ need only be an $[\ell,k,t]_q$ code.
\label{th:main}
\end{theorem}


\thm{th:main} can be strengthened further using the notion of an orthogonal
array \cite{HSS:ortho}.  Linear error-correcting codes are dual to linear
orthogonal arrays, and the proof actually uses orthogonal arrays rather than
codes.  In some cases non-linear orthogonal arrays are slightly better than
linear ones.  See Sections~\ref{s:proof} and \ref{s:other}.

The most effective case of \thm{th:main} is in the asymptotic limit $n \to
\infty$ with $t$ and $q$ fixed.  Recall that a function $f(n)$ is
\emph{quasilinear} if $f(n) = O((\log n)^\alpha n)$ for some $\alpha$.
Quasilinearity is also written $f(n) = \tO(n)$.  Say that a family $\{F\}$ of
cubature formulas is quasilinear (abbreviated ``QL'') if the points and weights
of each $F$ can be generated in quasilinear time in the length of the output.

\begin{theorem} Assume all variables as in \thm{th:main}.  Then $G$ can have
$O(\ell^\alpha)$ points (with the constant depending only on $q$), where
$$\alpha = t-1 - \bfloor{\frac{t-1}{q}}.$$
If each $G$ is centrally symmetric and $t$ is odd, then
$$\alpha = t-2 - \bfloor{\frac{t-2}{q}}.$$
Moreover, $G$ is quasilinear as $\ell \to \infty$, assuming precomputation of
each $F_i$.
\label{th:asymp}
\end{theorem}

If $\mu$ is an $m$-fold product with $m \propto n$, then we can take $\ell
\propto n$ in \thm{th:asymp}, so that $O(\ell^\alpha) = O(n^\alpha)$. In
comparison, the Tchakaloff upper bound is $O(n^t)$ points, or $O(n^{t-1})$ when
$t$ is odd and $\mu$ is centrally symmetric (\sec{s:other}).  Thus
\thm{th:asymp} is asymptotically better than Tchakaloff's theorem for all such 
product measures. Tchakaloff's theorem also does not guarantee equal weights.

Another comparison is with the cardinality of exact determination. A
$t$-cubature formula $F$ is \emph{overdetermined}, \emph{underdetermined}, or
\emph{exactly determined} if the parameters of its points provide fewer, more,
or the same number of degrees of freedom, respectively, as the constraints
imposed by integrating all polynomials of degree $t$.  The cardinality of exact
determination is $\Theta(n^{t-1})$ for general $\mu$ and $\Theta(n^{t-1})$ when
$t$ is odd and $\mu$ is centrally symmetric. Thus for product measures, the
formulas in \thm{th:asymp} are asymptotically exactly determined (up to a
constant factor that depends on $t$) when $q$ is large.  But when $q < t-1$, 
or $q < t-2$ in the odd and centrally symmetric case, they are asymptotically
overdetermined.

A third comparison is with the Stroud lower bound:  Any $t$-cubature formula in
$n$ dimensions, not necessarily interior or positive, requires
$\Omega(n^{\floor{t/2}})$ points.  \thm{th:asymp} achieves the Stroud bound (up
to a constant factor) when $q=2$.

A final comparison is with an interesting thinning construction of Novak and
Ritter for products of quadrature formulas \cite{NR:simple}. (It is similar to
an earlier construction due to Grundmann and M\"oller for the $n$-simplex
\cite{GM:combinatorial}.) They produce $t$-cubature formulas with $O(n^{\floor{t/2}})$
points, which is within a constant factor of the Stroud bound and better than
\thm{th:asymp} when $q>2$.  Crucially, their formulas are not positive,
although they can be made interior.  They also require that the factors of
$\mu$ be 1-dimensional.  The Novak-Ritter construction does generalize to
convolutions, as long as each factor formula has collinear points.

\thm{th:asymp} can be used to construct interesting cubature formulas for
several infinite sequences of regions and measures considered by Stroud
\cite[Ch 7,8]{Stroud:calc}:

\begin{theorem} For any $t$:
\begin{description}
\item[1.] The $n$-cube $C_n$ with uniform measure
has a QLEI $t$-cubature formula with $O(n^{\floor{t/2}})$ points.
\item[2.] The cubical shell $C_n - rC_n$ has a QLEI $t$-cubature formula with
$O(n^{\floor{t/2}+1})$ points.
\end{description}
For any odd $t \ge 3$:
\begin{description}
\item[1.] $\R^n$ with Gaussian weight function has a QLEI $t$-cubature formula
with $O(n^{t-2})$ points.
\item[2.] Any spherically symmetric measure on $\R^n$ has a QLPI $t$-cubature
formula with $O(n^{t-2})$ points. This includes the $n$-ball $B_n$, the
spherical shell $B_n - rB_n$, and the $(n-1)$-sphere $S^{n-1}$ with uniform
measure; and $\R^n$ with radial exponential weight function $\exp(-||\vx||_2)$.
\end{description}
For any $t \ge 2$:
\begin{description}
\item[1.] The $n$-simplex $\Delta_n$ has a QLPI $t$-cubature formula with
$O(n^{t-1})$ points.
\item[2.] The $n$-cross-polytope $C_n^*$ with uniform measure has a
QLPI $t$-cubature formula with $O(n^{\floor{3t/2}-1})$ points.
\end{description}
\label{th:stroud}
\end{theorem}

All cases of \thm{th:stroud} other than the cross polytope $C_n^*$ improve
the Tchakaloff bound.  On the other hand, the construction for the cube $C_n$
matches the Stroud bound up to a constant factor.  We admit that this
$t$-dependent factor is very generous when $t$ is large:  For each $t = 2s+1$,
it approaches $2\cdot s^s\cdot s!$ as $n \to \infty$ in the favorable case $n =
2^m$.  By contrast, the Novak-Ritter formulas use only $2^s$ more points than
the Stroud bound as $n \to \infty$.

\thm{th:stroud} partially solves a problem of Stroud \cite[p.
18]{Stroud:calc}:  Are there PI 5-cubature formulas for $C_n$, $B_n$, or
$\Delta_n$ with $O(n^2)$ or $O(n^3)$ points?  \thm{th:stroud} provides QLPI
5-cubature formulas with $O(n^2)$ points for $C_n$, $O(n^3)$ points for $B_n$,
and $O(n^4)$ points for $\Delta_n$.  In \sec{s:special}, we will establish a
special QLPI 5-cubature formula for $B_n$ with $O(n^2)$ points and QLPB and
QLPI 3-cubature formulas for $\Delta_n$ with $O(n)$ points.  Thus the only
remaining case of Stroud's question is the $n$-simplex in degree 4 or 5.

\begin{remark} The formula in \thm{th:stroud} for $S^{n-1}$ is technically a
QLPB formula if we take the definition of boundary in general topology. 
However, we take boundary in the sense of geometric topology, so that
\thm{th:stroud} is correct as stated.
\end{remark}

\acknowledgments

The author would like to thank Hermann K\"onig, Eric Rains and Hong Xiao for
useful discussions.  The author is also indebted to the late Arthur Stroud for
his excellent introduction to the cubature problem.

\section{Proofs}
\label{s:proof}

\begin{proof}[Proof of \thm{th:main}] First, identify an affinely
independent set of $q$ points in $\R^{q-1}$ with the finite field $\F_q$.
For each $1 \le i \le \ell$,
choose a linear map $\pi_i:\R^q \to \R^n$ that sends $\F_q$
to $F_i$, and define $\pi:\R^{(q-1)\ell} \to \R^n$ to be their direct sum:
$$\pi = \pi_1 \oplus \pi_2 \oplus \cdots \oplus \pi_\ell.$$
Because $F_1 * F_2 * \cdots * F_\ell$ is a $t$-design for the measure
$\mu$ on $\R^n$, the identity
$$\int P(\vx) d\mu = \frac1{q^\ell} \sum_{\vp \in \F_q^\ell}
    P(\pi(\vp))$$
holds for any polynomial $P$ of degree at most $t$ on $\R^n$.  Now suppose that
we thin the set $F = \pi(\F_q^\ell)$ to a set $G = \pi(A)$ for some set $A
\subset \F_q^\ell$.  Since $\pi$ is linear, if we want $G$ to be a $t$-cubature
formula for $\mu$ as $F$ is, it suffices that
\eq{e:A}{\frac1{q^\ell} \sum_{\vp \in \F_q^\ell} P(\vp)
    = \frac1{|A|} \sum_{\vp \in A} P(\vp)}
for any polynomial $P$ on $\R^{(q-1)\ell}$ of degree at most $t$. If $P$ is a
monomial, then as a function on $\F_q^\ell$ it depends on at most $t$
coordinates.  Conversely, any function on $\F_q^\ell$ it depends on at most $t$
coordinates is realized by a polynomial of degree at most $t$. It follows that
equation~\eqref{e:A} is equivalent to the statistical property that the
projection of $A$ onto any $t$ of the $I$ coordinates of $\F_q^\ell$ is
constant-to-1.  Such a set $A$ is called an \emph{orthogonal array} of 
strength $t$.

If $C$ is an $[\ell,k,t+1]_q$ code, then the dual space $C^*$ (in the sense of
linear algebra over $\F_q$) is a linear orthogonal array of strength $t$.
Since $C$ has dimension $k$, $C^*$ has dimension $\ell-k$ and therefore has
$q^{\ell-k}$ points.  Thus we can let $G = \pi(C^*)$.

The refinement when $t$ is odd and each $F_i$ is centrally symmetric is as
follows.  If $q$ is odd, we replace $\R^{q-1}$ by $\R^{(q-1)/2}$, and we
position $\F_q$ as a centrally symmetric set that does not lie in a
hyperplane.  (In other words, the points of $\F_q$ are the vertices of an
affinely regular cross polytope, plus the origin.)  We further demand that
negation in $\F_q$ coincides with negation in $\R^{(q-1)/2}$. Then any
centrally symmetric subset $A \subset \F_q^\ell$ is centrally symmetric in
$\R^{(q-1)\ell/2}$.  In this case both sides of \eqref{e:A} vanish when $P$ is
an odd polynomial. Thus $A$ need only be an orthogonal array of strength $t-1$.
In particular, this is so if $A = C^*$, because $C^*$ is a vector space over
$\F_q$ and vector spaces are centrally symmetric sets.

Finally if $t$ is odd, $q$ is even, and $C$ is a zero-sum code, then $C^*$
contains the vector $(1,1,\ldots,1)$ and is therefore invariant under addition
by this vector.  In this case we replace $\R^{q-1}$ in the general construction
by $\R^{q/2}$ and we realize $\F_q$ as a centrally symmetric set (the vertices
of a regular cross polytope). We further demand that adding 1 in $\F_q$
coincides with negation in $\R^{q/2}$.  Then once again $C^*$ is centrally
symmetric and need only be an orthogonal array of strength $t-1$. 
\end{proof}

The following lemma establishes \thm{th:asymp} as a corollary of \thm{th:main}.

\begin{lemma} Let $q$ be a prime power, let $m,t \in \Z_{\ge 0}$, and let
$$\alpha = t-1-\floor{\frac{t-1}{q}}.$$
Then there is a $[q^m,k,u]_q$ zero-sum code $C$ with
$$u \ge t+1 \qquad k \ge q^m-m\alpha-1.$$
\label{l:bch} \eatline \end{lemma}

The code in \lem{l:bch} is called an (extended, narrow-sense) \emph{BCH code}
\cite{SM:codes,CS:splag,BR:class,Hocquenghem:codes}.  We will use the duals of
BCH codes to thin cubature formulas. As it happens, the dual of a
BCH code of this type is another BCH of the same type.

\begin{proof} It is easier to define the dual code $C^*$ and show that it is an
orthogonal array.  Since it is a linear space, it suffices to show that every
coordinate projection $\pi_I:C^* \to \F_q^I$ with $|I| \le t$ is onto.  There
is an important $\F_q$-linear function
$$\Tr_q:\F_{q^m} \to \F_q$$
called the \emph{trace}.  (It is analogous to the taking the real part of a
complex number.)  First, we interpret $\F_q^{q^m}$ as the space of all
functions from $\F_{q^m}$ to $\F_q$. We define $C^*$ as the set of all
functions
$$f:\F_{q^m} \to \F_q \qquad f(x) = \Tr_q(P(x)),$$
where $P$ is a polynomial of degree at most $t-1$.  If $I \subseteq \F_{q^m}$
and $|I| \le t$, the polynomial $P$ can achieve any desired values on $I$ by
Lagrange interpolation.  Thus the distance of $C$ is at least $t+1$.

The space of polynomials of degree $t-1$ on $\F_{q^m}$ has $\F_{q^m}$-dimension
$t$, and therefore $\F_q$-dimension $mt$. But taking the trace reduces the
dimension in two ways. To give an explicit example, suppose that $q=2$, $t=3$, 
and $m$ is arbitrary.  Then $C^*$ is the set of all
$$f(x) = \Tr_2(ax^2+bx+c).$$
The apparent dimension of $C^*$ is $3m$.  But $f$ only depends on the trace of
$c$, so $c$ contributes $1$ rather than $m$ to the dimension of $C^*$.
Moreover, $\Tr_2(bx) = \Tr_2(b^2x^2)$, so the linear term can be removed from
$f$, with the conclusion that
$$\dim C^* \le m+1.$$
In general, the constant term of $P$ contributes $1$ to the dimension and the
other $t-1$ terms contribute $m$ each, except that $\floor{\frac{t-1}{q}}$
terms are superfluous by the Frobenius automorphism $x \mapsto x^q$.  Thus
$$\dim C^* \le m\alpha+1,$$
as desired.

Since constants are polynomials of degree $0$, $C^*$ contains constant
vectors.  Therefore $C$ is a zero-sum code.
\end{proof}

On the face of it, \lem{l:bch} only establishes \thm{th:asymp} when $\ell =
q^m$.  If $q^{m-1} < \ell < q^m$, we can project a BCH code from $\F_q^{q^m}$
to $\F_q^n$.  This preserves the $O(t^\alpha)$ bound at the expense of
worsening the constant factor.  If $\ell$ is not much more than $q^{m-1}$, we
can slightly improve the projected code with a projection that annihilates up
to $\alpha-1$ independent vectors in $C^*$. (See \thm{th:sphere5} for an
example.)

\begin{remark} The inequalities for $u$ and $k$ in \lem{l:bch}
become sharp as $m \to \infty$.
\end{remark}

\begin{proof}[Proof of \thm{th:stroud}] The simplest case to consider is with
uniform measure and $\R^n$ with Gaussian weight function. This fits
\thm{th:asymp} with $\ell=n$, provided that for each $t$, we find an EI
$t$-quadrature formula with Gaussian weight and with $q$ points for some prime
power $q$.  Since there is no bound on $q$, the Seymour-Zaslavsky theorem
\cite{SZ:averaging} establishes that such formulas exist.  One explicit method
begins with the PI Gaussian $(2t+1)$-quadrature formula with $t+1$ points. 
Viewed as a $t$-quadrature formula, the $t+1$ points can be freely perturbed. 
In particular they can be perturbed so that the weights become multiples of
$1/q$ for some large prime power $q$. The perturbation can be chosen to retain
central symmetry. On the other hand, since $q$ is large, $\floor{(t-2)/q} = 0$.
Thus \thm{th:asymp} produces formulas with $O(n^{t-2})$ points.

We will need the same construction for the orthant $\R_{\ge 0}^n$ with
exponential weight function $\exp(-||\vx||_1)$.  This measure does not have
central symmetry, and the end result is formulas with $O(n^{t-1})$ points,
again with $\ell=n$.

The next simplest case is the $n$-simplex $\Delta_n$. Recall that $\Delta_n$
has barycentric coordinates
$$x_0 + x_1 + \ldots + x_n = 1$$
which realize it as a subset of the orthant $\R_{\ge 0}^{n+1}$.
If $P(\vx)$ is a polynomial of
degree $t$ on $\Delta_n$, then it can be \emph{homogenized}:
it can be expressed as a homogeneous polynomial of degree $t$ by
attaching a factor of $\bigl(\sum_i x_i\bigr)^{t-s}$ to each
term of degree $s$.  In this case
$$\int_{\Delta_n} P(\vx) d\vx
    = \frac1{(n+t)!}\int_{\R_{\ge 0}^{n+1}} P(\vx) \exp(-||\vx||_1) d\vx.$$
Therefore we can project any non-exterior cubature formula for $\R_{\ge
0}^{n+1}$ radially onto $\Delta_n$ without loss of degree, although the
weights change.  (If the origin happens to be a cubature point, discard it.) 
In particular, we can project the cubature formulas provided by \thm{th:asymp}
as explained previously.  The formulas still have $O(n^{t-1})$ points,
although the weights are no longer equal.

The same argument works for the sphere $S^{n-1} \subset \R^n$ for centrally
symmetric formulas.  Every polynomial $P$ on $S^{n-1}$ can be expressed as
$P_S + P_A$, where $P_S$ is centrally symmetric and $P_A$ is centrally
antisymmetric.  The integral of $P_A$ vanishes, as does its sum with respect
to any centrally symmetric formula.  Meanwhile every term of $P_S$ has even
degree, so it can be expressed as a homogeneous polynomial on $\R^n$ using the
equation
$$x_1^2 + \ldots + x_n^2 = 1$$
for the unit sphere.  Then
$$\int_{S^{n-1}} P(\vx) d\Omega = \frac2{(\frac{n+t}2-1)!}
    \int_{\R^n} P(\vx) \exp(-||\vx||_2^2) d\vx,$$
where $\Omega$ is usual surface volume on $S^{n-1}$. Again, any centrally
symmetric cubature formula can be radially projected and the weights adjusted.

Formulas for the ball $B_n$ and the spherical shell $B_n - rB_n$ can be
derived from formulas for the sphere $S^{n-1}$ using radial separation of
variables \cite[Th 2.8]{Stroud:calc}. The result is a product formula where
the radial factor can be Gaussian quadrature.  The number of points in this
factor does not increase with dimension.

The cross-polytope $C_n^*$ is the union of $2^n$ simplices.  Thus we can obtain
formulas for $C_n^*$ by repeating formulas for $\Delta_n$.  In degree $t$, we
do not need all $2^n$ copies; instead we can repeat it in the pattern of the
BCH code over $\F_2$ defined by polynomials of degree $t-1$ over
$\F_{2^m}$. Such a code has $O(n^{\floor{t/2}})$ vectors and the formula for
$\Delta_n$ has $O(n^{t-1})$ points, so the total is $O(n^{\floor{3t/2}-1})$
points.

The $n$-cube $C_n = [-1,1]^n$ is in some ways the most interesting case.  Like
the  Gaussian case, it is a straight application of \thm{th:asymp} using an
equal-weight quadrature formula.  But in this case we will carefully
choose the quadrature formula on $[-1,1]$ to itself be a convolution
of $s = \floor{t/2}$ formulas with two points.  For example,
the Chebyshev 5-quadrature formula has points at
$$\pm \sqrt{\frac{5+\sqrt{5}}{30}} \pm \sqrt{\frac{5-\sqrt{5}}{30}}.$$
This is evidently a convolution, as is any centrally symmetric, equal-weight
formula with 4 points.  Elsewhere \cite{Kuperberg:moments} we show
that the $2^s$ points
$$\pm z_1 \pm z_2 \pm \cdots \pm z_s$$
form a Chebyshev-type $(2s+1)$-quadrature formula for $[-1,1]$ with constant
weight if and only if the $z_i$'s are the roots of the polynomial
$$Q(x) = x^s-\frac{x^{s-1}}{3} + \frac{x^{s-2}}{45} - \cdots +
    \frac{(-1)^s}{1\cdot3\cdot 15 \cdots (4^s-1)}.$$
We also show that all roots of $Q$ are real and that the resulting quadrature
formula is interior.  The $n$-fold product power of this formula is thus a
convolution of $sn$ pairs of points, so we can apply \thm{th:asymp} with $\ell
= sn$ and $q=2$.

Finally the $O(n^{\floor{t/2}})$ formula for the $n$-cube $C_n$ yields a
$O(n^{\floor{t/2}+1})$ formula for the cube surface $\partial C_n$ just by
repeating the formula for $C_{n-1}$ on each facet of $C_n$.  Then radial
separation of variables produces a product formula for the cubical shell $C_n
- rC_n$ which also has $O(n^{\floor{t/2}+1})$ points.
\end{proof}

\section{Special constructions and examples}
\label{s:special}

In this section we will consider some examples and special constructions with
concern for constant factors. For this purpose, we spell out more precisely the
notion of an orthogonal array.  Let $A$ be a finite set.  If a subset $X
\subset A^n$ has the property that its projection $X \to A^I$ is a
constant-to-1 map for every $|I| \le t$, then $T$ is an \emph{orthogonal array}
of strength $t$, or an $\OA(|T|,n,|A|,t)$ \cite{HSS:ortho}.  If $A = \F_q$
and $X = C^*$ is the dual of an $[n,k,t]_q$ code, then
$X$ is an $\OA(q^{n-k},n,q,t-1)$.  We will also say that $X$
is an $[n,n-k,t^*]_q$ to refer to its linear structure and indicate
its dual distance.

If $|S| = q$ is a prime power and $t$ is fixed, then BCH codes are the
best presently known $\F_q$-linear orthogonal arrays in the limit $n \to
\infty$.  But a few non-linear arrays are slightly better.  

A \emph{Hadamard matrix} of order $n$ is an $n \times n$ matrix with entries
$\pm 1$ and with orthogonal rows (and therefore orthogonal columns as well). 
It is easy to show that a Hadamard matrix is equivalent to an $\OA(2n,n,2,3)$. 
A $[2^m,m+1,4^*]_2$ BCH code, which is also called a first-order
Reed-Muller code, yields a Hadamard matrix of order $2^m$.  But there are also
Hadamard matrices for other values of $n$, for example when $4|n$ and $n-1$ is
prime.  The Hadamard conjecture asserts that there is a Hadamard matrix of
every order $n$ divisible by $4$.

For any even $m \ge 4$, there is a Kerdock code which is a non-linear
$\OA(2^{2m},2^m,2,5)$. It has $\frac12$ as many points as the corresponding
$[2^m,2m+1,6^*]_2$ BCH code
\cite{HSS:ortho,HKCSS:linearity,Kerdock:codes}. For any even $m \ge 6$, there
is a Delsarte-Goethals code which is a non-linear $\OA(2^{3m-2},2^m,2,7)$
\cite{DG:alternating}.  It has $\frac14$ as many points as the corresponding
$[2^m,3m+1,8^*]_2$ BCH code.

The following result comes from thinning some cubature formulas of Stroud, some
of whose points have a product structure.

\begin{theorem} Let $n \ge 6$ and let
$$k = \begin{cases}
    4m & 2^{2m-1} < n \le 2^{2m} \\
    4m+2 & 2^{2m} < n \le 2^{2m}+2^m \\
    4m+3 & 2^{2m}+2^m < n \le 2^{2m+1}
    \end{cases}$$
Then the sphere $S^{n-1}$, $\R^n$ with Gaussian measure, and the ball $B_n$
admit QLPI 5-cubature formulas with $2^k+2n$ points.
\label{th:sphere5}
\end{theorem}

\begin{proof}
The formulas $S_n$:5-3, $U_n$:5-2 and $E_n^{r^2}$:5-3 listed in Stroud
\cite[pp. 270,294,317]{Stroud:calc} have $2^n+2n$ points with $2^n$ of them
lying on the vertices of a cube.  These $2^n$ points can be thinned to either
the $[2^{2m+1},4m+3,6^*]_2$ BCH code, or the Kerdock
$\OA(2^{4m},2^{2m},2,5)$, and then projected down to $n$ dimensions.

If $2^{2m} < n \le 2^{2m}+2^m$, then the $[2^{2m+1},4m+3,6^*]_2$ BCH
code can be reduced by half by carefully choosing the projection. The code has
a vector of weight $2^{2m}-2^m$, so when $n$ is only slightly larger than
$2^{2m}$, we can choose a projection that annihilates this vector.

In each of the three cases, the result is a formula with $2^k+2n$ points.
\end{proof}

Actually, \thm{th:sphere5} is not quite optimal, because it uses a convenient
set of good distance-6 linear codes and non-linear strength-5 orthogonal arrays
rather than the best ones presently known.  A complicated map of the best
presently known linear codes over $\F_2$ of length $n \le 256$ is provided by
the ``best codes" functions in Magma \cite{Magma}.  Undoubtedly this map
could be augmented by non-linear orthogonal arrays, but we know of no effort to
do so.  When $n$ is a power of 2, Kerdock and BCH codes are the best presently
known choices.

Victoir \cite{Victoir:cubature} also established \thm{th:sphere5} (with
BCH codes). If $n = 2^m$ and Stroud's formulas for $S^{n-1}$ is
thinned using a BCH code, it then has equal weights and is therefore
a $5$-design. Interestingly, in this case it has a transitive symmetry group
and was previously found by Calderbank, Hardin, Rains, Shor, and Sloane
\cite{CHRSS:grass}.  Similar constructions were found by K\"onig
\cite{Konig:cubature}, by Sidelnikov \cite{Sidelnikov:7des}, and by Schechtman,
interpreting work of Hajela \cite{Hajela:thin}.

We can obtain a good 3-cubature formula for the cube $C_n$ by a straightforward
application of \thm{th:asymp} using the 2-point Gaussian quadrature formula for
the interval $[-1,1]$.  Thinning the product formula using a BCH code
yields a $2^{j+1}$-point formula when $2^{j-1} < n \le 2^j$. When $n=2^j$, or
more generally whenever there is a Hadamard matrix of order $n$, the product
formula can be thinned to  the $2n$ vertices of a certain regular
cross-polytope inside $C_n$. A formula due to Stroud ($C_n$:3-1 \cite[p.
230]{Stroud:calc}) also uses the vertices of a regular cross-polytope, but not
the same one.

We can obtain a 3-cubature formula with $O(n)$ points for $\Delta_{n-1}$ with a
similar construction.  Using known Hadamard matrices, the formula has $3n+o(n)$
points; if the Hadamard conjecture holds, it has between $3n-1$ and $3n+5$
points.  First, the positive ray $\R_{\ge 0}$ with exponential weight has a
equal-weight 2-quadrature formula with points at $0$ and $2$.  If we apply
\thm{th:main} to this formula and a Hadamard matrix of order $n$, the result is
a $2n$-point formula $F$ on $\R_{\ge 0}^n$ which also has degree 2.  However,
if our interest is integration on $\Delta_{n-1}$, we need only consider
homogeneous polynomials on $\R_{\ge 0}^n$. The formula $F$ correctly integrates
every degree 3 monomial other than $x_i^3$.  We can fix $F$ for these
monomials, without changing its sum for $x_i^2x_j$ or $x_ix_jx_k$, by adding a
point at $(1,0,0,\ldots)$ (and permutations) with weight $2$.

The projected formula on $\Delta_{n-1}$ consists of these points and
weights in barycentric coordinates:
$$\begin{array}{cl}
(1,0,0,\ldots,0)_S & \frac{2}{n(n+1)(n+2)} \\
(\frac2n,\frac2n,\ldots,\frac2n,0,0,\ldots,0)_H & \frac{n}{2(n+1)(n+2)} \\
(\frac1n,\frac1n,\ldots,\frac1n) & \frac{4n}{(n+1)(n+2)}
\end{array}.$$
The subscript ``S'' denotes full symmetrization, as in Stroud's notation.  The
subscript ``H'' denotes symmetrization in the pattern of a Hadamard design.
(See \sec{s:other}.)  This produces a formula with $3n-1$ points provided that
there exists a Hadamard matrix of order $n$. When there is none, we can use a
Hadamard matrix of order $\ell > n$.  The formula on $\Delta_{\ell-1}$ with
$3\ell-1$ points can be projected onto $\Delta_{n-1}$, as in the proof of
\thm{th:stroud}. We can take $\ell = n + o(n)$ by letting $\ell-1$ be the first
prime after $n$ which is 3 mod 4.  If the Hadamard conjecture holds, we can
take $\ell = 4\ceil{n/4}$.

Stroud asked for a practical, PI 5-cubature formula for $C_{100}$.  Following
\thm{th:stroud}, we can find one by thinning the product formula coming from
the 4-point Chebyshev quadrature on $[-1,1]$. This product formula is the
convolution of 200 pairs of points, so we can thin it using the Kerdock
$\OA(2^{16},2^8,2,5)$, projected to $200$ dimensions. The cubature formula
therefore has $2^{16} = 65536$ points, which would have been fairly practical
even in 1971 when Stroud asked the question.  (The Kerdock code used here was
discovered shortly afterward \cite{Kerdock:codes}, but the BCH codes was known
in 1959 \cite{BR:class,Hocquenghem:codes}.)

Victoir \cite{Victoir:cubature} found another thinning of the same Chebyshev
product formula with $4^{12} = 16777216$ points, which the author tied in the
first version of this paper. 

Note that the Chebyshev-Kerdock 5-cubature formula for $C_{100}$ is
overdetermined. The threshold of exact determination for centrally symmetric
5-cubature formulas on $C_{100}$ is 87651 points.  Meanwhile the centrally
symmetric Tchakaloff bound is 8852652 points, while the Stroud lower bound is
5050 points.

Finally Sch\"urer \cite{Schurer:comparison} compared the numerical accuracy of
various cubature and quasi-Monte-Carlo methods for the integration of various
test functions defined on $C_n$ with $2 \le n \le 100$.  He assumed a more
modern limit of $2^{25}$ evaluations of the integrand. For much of this test
regime we can suggest the following cubature formulas:  Start with the power of
the convolutional 7-quadrature formula \cite{Kuperberg:moments} for
$[-1,1]$, whose points are approximately at
$$\pm .500128 \pm .243941 \pm .153942.$$
Then thin the $n$-fold product power of this formula using a Delsarte-Goethals
code.  The result is an EI 7-cubature formula with at most $2^{23}$ points up
to dimension $\floor{256/3} = 85$.

\section{Other comments}
\label{s:other}

Victoir \cite{Victoir:cubature} proposes thinning symmetric cubature formulas
rather than product or convolution formulas.  The enabling result of symmetric
cubature formulas is Sobolev's theorem:  If a linear action of a finite group
$G$ preserves $\mu$, then a cubature formula consisting of orbits of $G$ need
only be checked for $G$-invariant polynomials. Victoir extends Sobolev's
theorem with a $G$-invariant generalization of Tchakaloff's theorem:  A PI
cubature formula only needs as many orbits as the dimension of $\R[\vx]^G_{\le
t}$, the space of $G$-invariant polynomials of degree at most $t$.   One
important special case is when $G$ is the 2-element central symmetry group.  If
$\mu$ is a measure on $\R^n$ with central symmetry and $t$ is odd, the bound
from this version of Tchakaloff's theorem is $O(n^{t-1})$ points.

Even if a cubature formula $F$ uses very few orbits of $G$, some of the orbits
might be very large.  Victoir proposes thinning each large orbit separately. He
notes that this can be done using linear programming, among other methods;
linear programming on a set of $G$-orbits should be much easier than general
numerical methods to find positive cubature formulas for $\mu$. If $G =
(\Z/2)^n$ is the group of independent sign changes of all $n$ coordinates, then
an orbit of $G$ is a Cartesian power can be identified with $\F_2^k$ for some
$k \le n$. In this case Victoir found the constructions of \thm{th:main} and
\thm{th:asymp}.  (In the case of 5-cubature on $C_n$, he found a special 
construction with $O(n^3)$ points with elements of both \thm{th:main} and the
$n$-cube case of \thm{th:stroud}.)

If $G$ is the group of coordinate permutations, then an orbit whose points have
two distinct coordinates can be identified with the set of $k$-subsets of an
$n$-set. A geometric $t$-design $T$ within this orbit is also a traditional
combinatorial $t$-design, or an $(n,k,t)-\lambda$ design.  Namely, $T$ is a
collection of blocks of size $k$ in a set of $n$ such that each $t$-subset is
contained in exactly $\lambda$ blocks.  In particular, an
$(n,\frac{n}2,3)-\frac{n}4$ design is called a \emph{Hadamard design}, because
it comes from the rows of a Hadamard matrix.

These constructions motivate the notion of a weighted orthogonal array.  We
define it as a finite set $A$ and a measure $\mu$ on $A^n$ that projects to
uniform measure on each $A^I$ with $|I| \le t$.  More generally, $\mu$ might
project to $\sigma^n$ for some reference measure $\sigma$ on $S$.  Such arrays
could improve of \thm{th:main}; the factor formulas would not need to have
equal weights.

Finally, cubature formulas coming from \thm{th:main} could be viewed as
quasi-Monte-Carlo methods.  They are similar to some constructions of
$(t,m,s)$-nets, which are quasi-Monte-Carlo methods first defined and largely
developed by Niederreiter \cite{Niederreiter:tms}.  Nonetheless PI cubature
formulas and discrepancy-based quasi-Monte-Carlo methods are thought to have
complementary advantages \cite{Schurer:comparison}.  We believe that
the improved asymptotics presented here could change the standing of cubature
among numerical methods for integration.


\providecommand{\bysame}{\leavevmode\hbox to3em{\hrulefill}\thinspace}

\end{document}